\documentclass{iau}
\usepackage{graphicx}
\usepackage{subcaption}
\usepackage{amsmath, amsbsy, amsfonts}
\usepackage{natbib}

\def\Lambdavec{\boldsymbol{\Lambda}}
\def\lambdavec{\boldsymbol{\lambda}}
\def\xivec{\boldsymbol{\xi}}
\def\etavec{\boldsymbol{\eta}}

\title{On the relativistic Lagrange-Laplace secular dynamics for extrasolar systems}

\author[M. Sansottera, L. Grassi and A. Giorgilli]{M. Sansottera$^1$, L. Grassi$^1$, and A. Giorgilli$^1$}

\affiliation{$^1$Dipartimento di Matematica,\\
Universit\`a degli Studi di Milano, \\
20133 --- Milano, Italia \\
\begin{tabular}{l l}email: &{\tt marco.sansottera@unimi.it} \\
&{\tt lorenzo.grassi@studenti.unimi.it} \\
&{\tt antonio.giorgilli@unimi.it}
\end{tabular}
}

\pubyear{2015}
\volume{S310}  
\pagerange{74--77}
\jname{Complex Planetary Systems}
\editors{Z. Kne\v{z}evi\'c and A. Lema\^{\i}tre eds.}
\doi{10.1017/S174392131400787X}
\begin{document}

\maketitle

\begin{abstract}
We study the secular dynamics of extrasolar planetary systems by
extending the Lagrange-Laplace theory to high order and by including the
relativistic effects.  We investigate the long-term evolution of the
planetary eccentricities via normal form and we find an excellent
agreement with direct numerical integrations.  Finally we set up a
simple analytic criterion that allows to evaluate the impact of the
relativistic effects in the long-time evolution.

\end{abstract}

\firstsection 
\section{Introduction}
The study of the secular evolution of planetary systems is a long
standing and challenging problem.  The discoveries of hundreds of
extrasolar planetary systems raised many interesting problems
concerning their long-term evolution.  In the present paper, we study
the secular dynamics of two non-resonant coplanar planets in an
extrasolar system.  We extend the Lagrange-Laplace theory to high
order and include the main relativistic effects.

The study of extrasolar system raised two particularly relevant
problems, namely: (i)~most exoplanets have highly eccentric orbits, in
contrast with the almost circular orbits of the Solar System;
(ii)~there are many giant planets orbiting at a low distance from the
central star, with periods of a few months or even a few days.  In the
latter case relativistic effects could have a significant impact and
should be taken into account.  The General Theory of Relativity,
despite having been widely used in astrophysics, is not commonly
adopted in the study of planetary system dynamics.

The generalization of the Lagrange-Laplace secular theory to high
order in the eccentricities has been exploited so as to obtain an
analytic model that gives an accurate description of the behavior of
planetary systems, up to surprisingly high eccentricities (see,
\cite{LibHen-2005,LibHen-2006}). The results appear to be quite good
for systems which are not close to a mean-motion resonance. In
\cite{LibSan-2013} the secular theory has also been extended to order
two in the masses, by using a first-order approximation of an elliptic
lower dimensional torus in place of the usual circular approximation.
In particular this allows to deal with systems close to a mean-motion
resonance. The relevance of the relativistic corrections and tidal
effects on the long-term evolution of extrasolar planetary systems
has been studied, e.g., in~\cite{AdaLau-2006} and \cite{MigGoz-2008}.

On the other hand, the application of Kolmogorov and Nekhoroshev
theorems, allowed to make substantial progress for the problem of
stability of the Solar System.  Indeed, in recent years, the estimates
for the applicability of both theorems to realistic models of some
part of the Solar System have been improved by some authors
(e.g.,~\cite{Robutel-1995}, \cite{CelChi-2005}, \cite{LocGio-2007},
\cite{GioLocSan-2009, GioLocSan-2014} and \cite{SanLocGio-2011b,
  SanLocGio-2011a}).

In the present paper we exploit the idea of extending the
Lagrange-Laplace theory, already used in the above-cited papers, to
the case of high eccentricities.  The technical tool is the
construction of a suitable normal form which allows us to investigate
the long-time evolution of the planetary eccentricities.  In this
contribution we neglect the tidal effects, although we know that for
many system they can be relevant.  We decided to just consider the
relativistic correction in order to keep the discussion at a simple
level and to show that the extension of the Lagrange-Laplace theory,
including relativistic effects, produces accurate results. We plan to
further investigate the problem in a forthcoming work.

\section{Classical expansion of the Hamiltonian}
We consider a system of three coplanar point bodies, mutually
interacting according to Newton's gravitational law: a central star
$P_0$ of mass $m_0$ and two planets $P_1$ and $P_2$ of mass $m_1$ and
$m_2$ and semi-major axis $a_1$ and $a_2$, respectively.

We refer to \cite{LibSan-2013} for a detailed exposition concerning
the expansion of the Hamiltonian in the Poincar\'e canonical variables,
that reads
\begin{equation}
H(\Lambdavec, \lambdavec, \xivec, \etavec)=
H_0(\Lambdavec)+\varepsilon H_1(\Lambdavec, \lambdavec, \xivec
,\etavec)\ ,
\label{eq:H_iniz_poinc}
\end{equation}
where $H_0$ is the Keplerian part and $\varepsilon H_1$ the perturbation
due to the mutual attraction between the planets.  Using the standard
notation, we will refer to $(\Lambdavec,\lambdavec)$ as the {\it fast
  variables} and to $(\xivec,\etavec)$ as the {\it secular variables}.

\section{Relativistic corrections}
Starting from the Hamiltonian of the Newton model, we add the
relativistic corrections due to the mutual interaction between the
star and each of the two planets. That is, we consider the correction
included in the relativistic Hamiltonian of the problem of two-body in
heliocentric coordinates $(\mathbf r, \mathbf p)$. The relativistic Hamiltonian takes the form
\begin{equation}
H=H_0+\varepsilon H_1 + \frac{1}{c^2}H_2,
\end{equation}
with $H_0$ and $\varepsilon H_1$ as in the Newtonian model, while $\frac{1}{c^2}H_2$ is
\begin{equation}
\frac{1}{c^2} H_2 = \frac{1}{c^2} \sum_{i=1}^2 \left[ -\frac{\gamma_{1, i}}{\mu_i^3} (\mathbf P_i \cdot \mathbf P_i)^2 - \frac{\gamma_{2, i}}{\mu_i} \frac{\mathbf P_i \cdot \mathbf P_i}{\parallel \mathbf r_i \parallel} - \frac{\gamma_{3, i}}{\mu_i} \frac{(\mathbf r_i \cdot \mathbf P_i)^2}{\parallel \mathbf r_i \parallel ^3 } +\gamma_{4, i} \mu_i \frac{1}{\parallel \mathbf r_i \parallel^2 } \right],
\end{equation}
with
\begin{equation}
\begin{split}
& \mu_i=\frac{m_0 m_i}{m_0 + m_i}, \qquad \beta_i=\mathcal G (m_0 + m_i), \qquad \upsilon_i=\frac{m_0 m_i}{(m_0 + m_i)^2}\ ,\\
& \gamma_{1, i}= \frac{1-3 \upsilon_i}{8}, \qquad \gamma_{2, i}=\frac{\beta_i(3+\upsilon_i)}{2}, \qquad \gamma_{3, i}=\frac{\beta_i  \upsilon_i}{2}, \qquad \gamma_{4, i}=\frac{\beta_i^2}{2}\ , 
\end{split}
\end{equation}
and $\mathbf P_i=\mathbf p_i + \frac{\mu_i}{m_0}\mathbf p_{3-i} + \mathcal{O}(c^{-2})$ for $i=1,2$. 


\begin{figure}
\begin{center}
  \begin{subfigure}{.4\textwidth}
    \includegraphics[width=\textwidth]{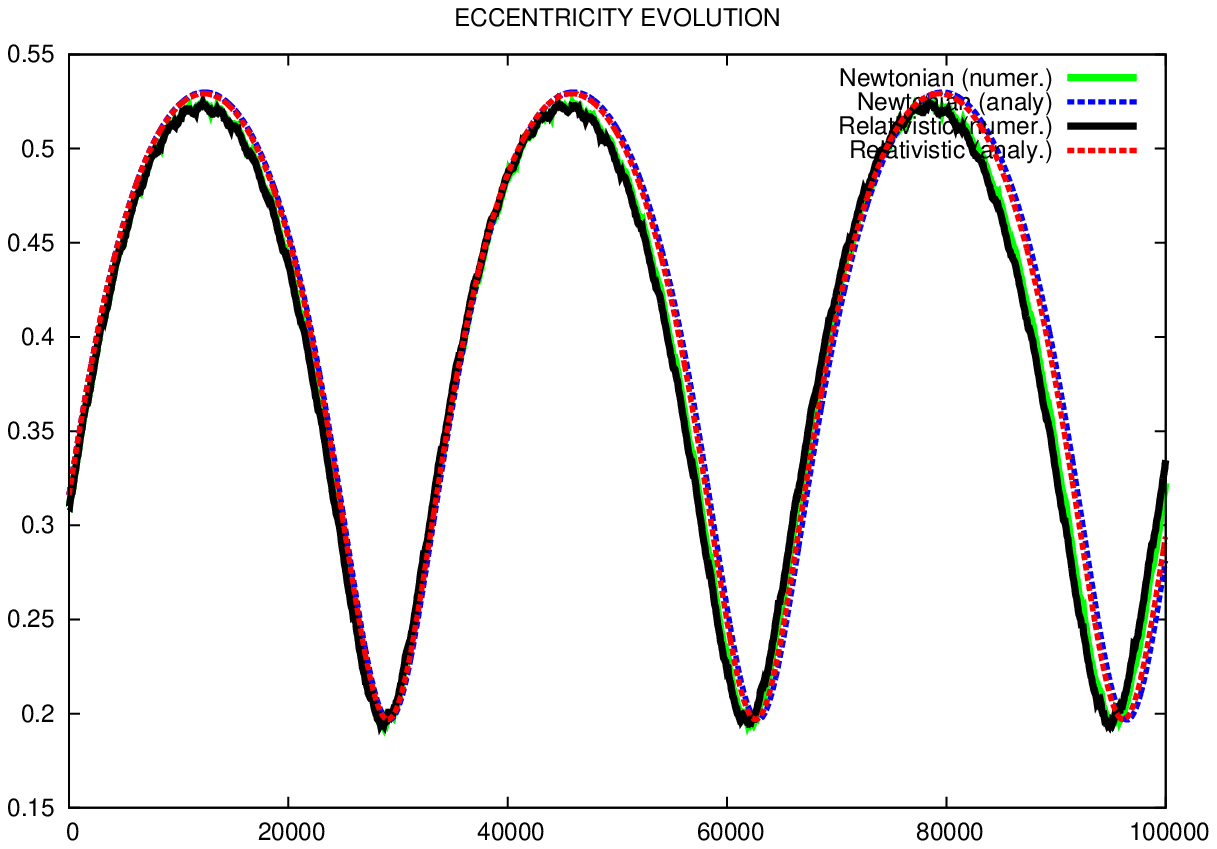}
    \caption{HD~169830, inner planet}
    \label{fig-11}
  \end{subfigure}
  \begin{subfigure}{.4\textwidth}
    \includegraphics[width=\textwidth]{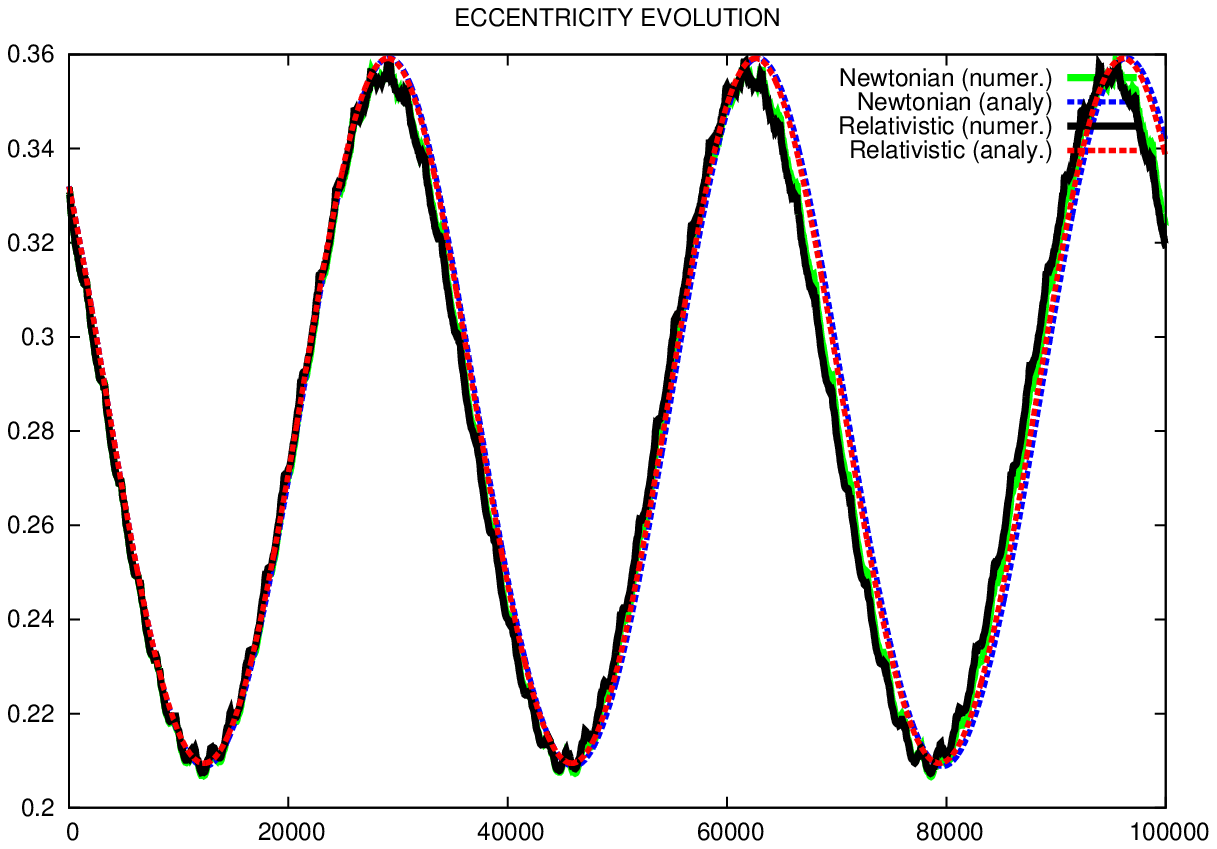}
    \caption{HD~169830, outer planet}
    \label{fig-12}
  \end{subfigure}

  \begin{subfigure}{.4\textwidth}
    \includegraphics[width=\textwidth]{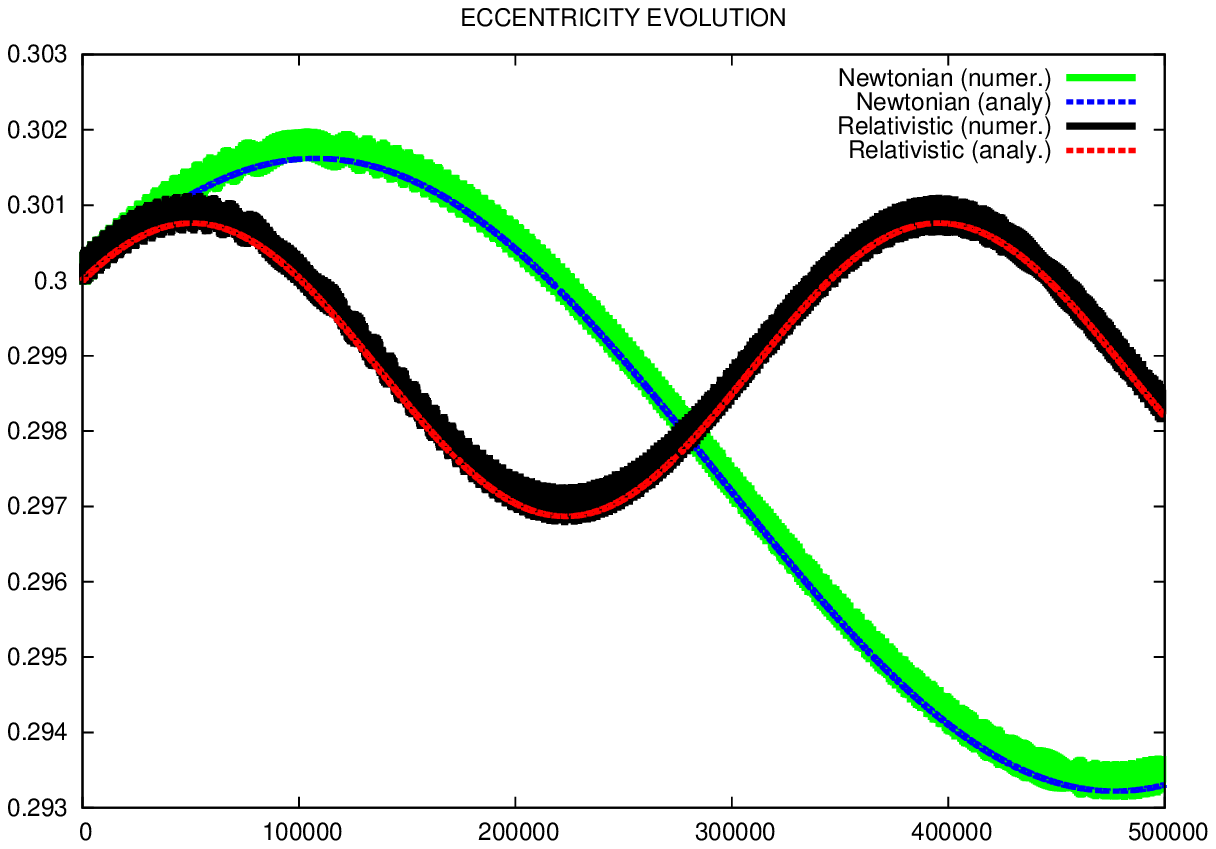}
    \caption{HD~11964, inner planet}
    \label{fig-21}
  \end{subfigure}
  \begin{subfigure}{.4\textwidth}
    \includegraphics[width=\textwidth]{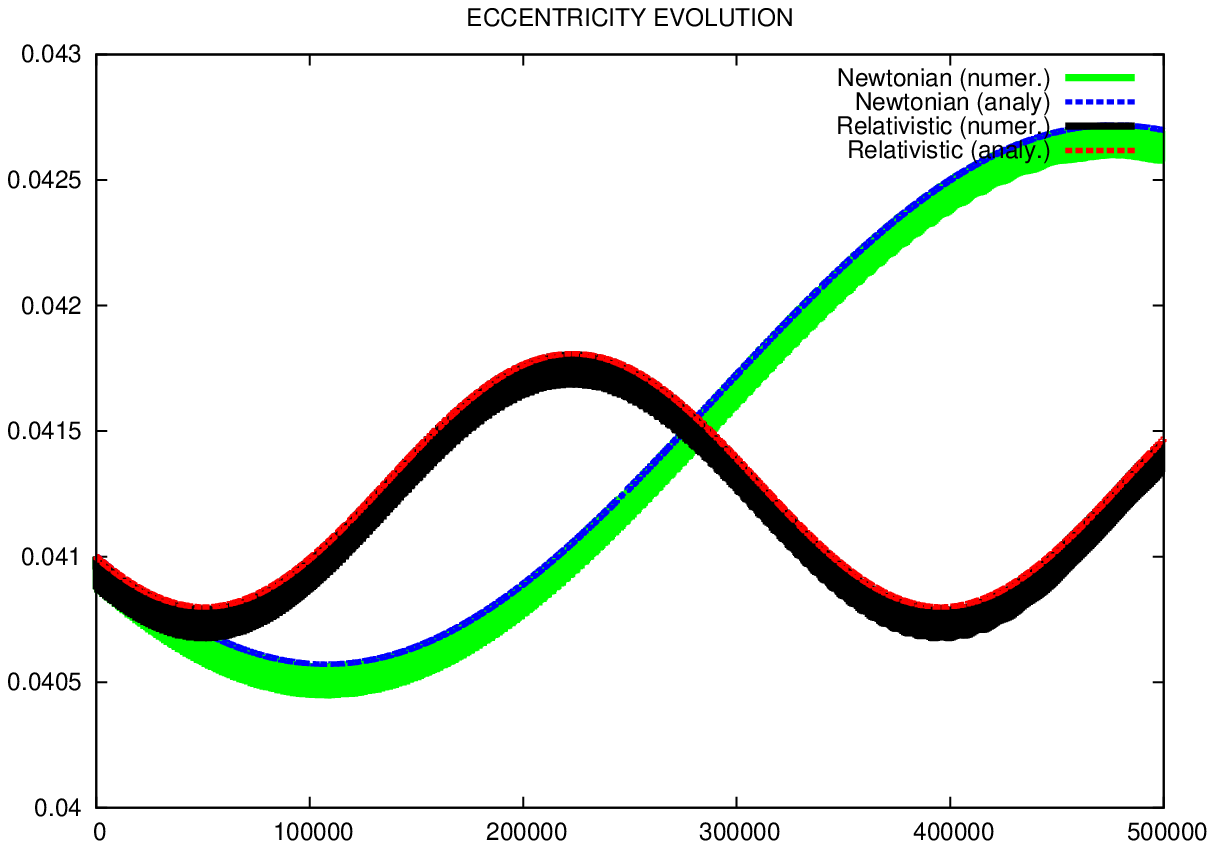}
    \caption{HD~11964, outer planet}
    \label{fig-22}
  \end{subfigure}
  \caption{Long-term evolution of the eccentricities for the HD~169830
    and HD~11964 planetary systems.  Comparison of the results
    obtained via direct numerical integration (green-black) against
    normal form (blue-red), for the Newtonian approximation
    (green-blue) and the model including the relativistic effects
    (black-red).}
  \label{fig1}
\end{center}
\end{figure}

\section{Long-term evolution}
As we are interested in the long-term dynamics, we remove the
dependency of the Hamiltonian from the fast angles.  The classical
approach consists in replacing the Hamiltonian with its average, the
so-called approximation at order one in the masses.  We replace this
procedure by a Kolmogorov-like step, that allows us to include in the
secular model the effects of the main near-resonances effects
(see~\cite{LibSan-2013} for a detailed exposition).  This is the
secular Hamiltonian at order two in the masses.

After averaging, the secular Hamiltonian has two degrees of freedom
and its quadratic part differs from the one considered in the
Lagrange-Laplace theory by relativistic corrections and contributions
of order two in the masses, which however are small.

As in \cite{LibSan-2013}, we introduce
the action-angle variables via normal form.  In the normalized
coordinates, the equations of motion take a simple form that can be
analytically integrated.  We validate the results by comparing the
analytic integration with the direct numerical integration of the
full three-body system.

In Figure~\ref{fig-11}--\ref{fig-12} we report the results for the
HD~169830 system.  In this case the relativistic effects are
negligible and the Newtonian approximation allows to accurately
describe the long-term evolution for a time interval of $10^5$ years.
Instead, for the HD~11964, the relativistic corrections play a major
role, as it is clearly shown in Figure~\ref{fig-21}--\ref{fig-22}.  In
this case the calculations cover $5\times10^5$ years. In all cases,
the evolutions via normal and via numerical integration are in
excellent agreement.  We emphasize that the use of normal form
provides us with a natural criterion for deciding whether or not the
relativistic corrections are relevant.  Indeed the difference is seen
in the precession frequencies: if the relativistic corrections are
relevant then so is the difference, as the figures clearly show.

\section{Relevance of the relativistic corrections}
In order to evaluate the impact of the relativistic corrections, we
look at the quadratic parts of the secular Hamiltonians, namely
\begin{equation*}
H_q^\text{(New)}(\boldsymbol \eta, \boldsymbol \xi) = \boldsymbol \eta
\cdot A \boldsymbol \eta + \boldsymbol \xi \cdot
A \boldsymbol \xi \quad\hbox{and}\quad H_q^\text{(Rel)}(\boldsymbol \eta,
\boldsymbol \xi) = \boldsymbol \eta \cdot B \boldsymbol
\eta + \boldsymbol \xi \cdot B \boldsymbol \xi\ ,
\end{equation*}
where $A$ and $B$ are real symmetric $2 \times 2$ with
$$
B=A - \frac{3}{2}\frac{ \mathcal G^{3/2}}{c^2}
\begin{bmatrix}
  \frac{(m_0+m_1)^{3/2}}{a_1^{5/2}} & 0\\
  0&\frac{(m_0+m_2)^{3/2}}{a_2^{5/2}}
\end{bmatrix}\ .
$$
Clearly, the relativistic effects are more important if
\begin{equation}
A_{ii}\sim -\frac{3}{2}\frac{ \mathcal G^{3/2}}{c^2} \frac{(m_0+m_i)^{3/2}}{a_i^{5/2}},
\qquad\hbox{i.e. if}\quad
\Pi_i \equiv \frac{4 \mathcal G a_2^3 m_0 (m_0+m_i)}{c^2 a_i^2 a_1^2 m_{3-i}} \sim 1.
\end{equation}
In the following table we report the dimensionless quantities $\Pi_i$ for the extrasolar systems considered above.

{
\renewcommand{\arraystretch}{1.2}
\setlength{\tabcolsep}{10pt}
\centerline{
\begin{tabular}{l c c}
HD~169830 & $\Pi_1$: 0.0021779 & $\Pi_2$: 0.0001547\cr
HD~11964  & $\Pi_1$: 0.9651708 & $\Pi_2$: 0.0399271\cr
\end{tabular}}
\vskip1pt
}
\noindent
We observed that for the majority of the extrasolar systems taken into
consideration, the relevance of relativistic corrections may be
inferred from the difference between the matrices.  This
provides us with a rough criterion based on the first order
approximation.  Normal form provides a more refined criterion.

\end{document}